\documentclass[12pt]{amsart}
\usepackage{amsmath, amsthm, amssymb, amsfonts}
\usepackage[english]{babel}
\usepackage{geometry}
\usepackage{xcolor}
\usepackage{enumitem}
\usepackage{bbm}
\usepackage{indentfirst}
\usepackage{amsaddr}
\usepackage{pifont}
\usepackage{xcolor}
\usepackage{hyperref}

\theoremstyle{definition}
\newtheorem{lemma}{Lemma}
\newtheorem{theorem}[lemma]{Theorem}
\newtheorem{definition}[lemma]{Definition}
\newtheorem{proposition}[lemma]{Proposition}
\newtheorem{corollary}[lemma]{Corollary}

\newtheorem{question}[lemma]{Question}
\DeclareMathOperator{\andd}{and}

\newcommand{\rest}{\upharpoonright} 

\def\cA{{\mathcal{A}}} \def\cB{{\mathcal{B}}}    \def\cF{{\mathcal{F}}}   \def\cI{{\mathcal{I}}}       \def\cP{{\mathcal{P}}}          \def\cZ{{\mathcal{Z}}}

\def\w{{\omega}}

\title{New examples of MAD families with pseudocompact hyperspaces}

\author[C. Corral and V. O. Rodrigues]{Cesar Corral and Vinicius de O. Rodrigues}

\address[C. Corral, V. O. Rodrigues]{York University}
\email{cicorral@yorku.ca and vor@yorku.ca}
\makeatletter
\renewcommand{\email}[2][]{%
	\ifx\emails\@empty\relax\else{\g@addto@macro\emails{,\space}}\fi%
	\@ifnotempty{#1}{\g@addto@macro\emails{\textrm{(#1)}\space}}%
	\g@addto@macro\emails{#2}%
}
\makeatother
\begin{document}
	\maketitle

\begin{abstract}
    We show that both $\mathfrak{ap=c}$ and $\diamondsuit(\mathfrak{b})$ imply the existence of MAD families with  pseudocompact Vietoris hyperspace, substantially expanding the list of models where their existence is known. We also discuss some properties on the structure of fin-intersecting MAD families by providing some new examples of non-fin-intersecting MAD families with special properties. We show that the Baire number of $\omega^*$ is strictly larger than $\mathfrak c$ if and only if for every MAD family $\mathcal A$, $Psi(\mathcal A)$ and its Vietoris hyperspace are $p$-pseudocompact for some free ultrafilter $p$.
\end{abstract}

\section{Introduction}

Topologies on spaces of subsets of a given topological space have been studied since the beginning of the last century. An early and important such topology is the Vietoris Hyperspace:
Let $X$ be a $T_1$ Topological space, by $\exp(X)$ we denote the set of all nonempty closed subsets of $X$.
Given a subset $U$ of $X$, we let 
$$U^+=\{F\in \exp(X):F\subseteq U\}\textnormal{ and }U^-=\{F\in \exp(X):F\cap U\neq \emptyset\}.$$
The \textit{Vietoris topology} on $\exp(X)$ is the topology generated by the subbase $\{U^+, U^-: U\subseteq X\text{ is open}\}$.
The \textit{Vietoris hyperspace} of $X$ is the set $\exp(X)$ endowed with the Vietoris topology. We denote by $\langle U_0,\ldots,U_n\rangle$ the basic open set
$$\bigg(\bigcup_{i\leq n}U_i\bigg)^+\cap \bigcap_{i\leq n}U_i^-$$
	
The Vietoris Topology was introduced by L. Vietoris \cite{vietoris1922bereiche} and coincides with the topology generated by Hausdorff's metric in case $X$ is a compact metric space. A classical survey on the basics of the Vietoris topology is \cite{michael1951topologies}.

Vietoris also proved that $X$ is compact if and only if its Vietoris hyperspace is compact \cite{vietoris1922bereiche}, raising the natural question of whether there are similar results for generalizations of compactness. Generalizations of this result were investigated by J. Ginsburg, who proved some results in this direction regarding pseudocompactness, countably compactness, $p$-compactness and $p$-pseudocompactness (where $p$ is some fixed free ultrafilter) \cite{ginsburg1975some}. He asked whether there is a relation between the pseudocompactness of $X^\omega$ and of $\exp(X)$ and whether there is a characterization for the spaces $X$ which have pseudocompact hyperspace.

A space $X$ is said to be pseudocompact if every continuous function from $X$ into $\mathbb R$ is bounded. Moreover, if $p \in \omega^*$ (the set of all free ultrafilters in $\omega)$, $X$ is said to be $p$-pseudocompact if every sequence of nonempty open sets has a $p$-limit point, that is, for every family of nonempty open sets $(U_n: n \in \omega)$ of $X$, there exists $x \in X$ such that $\{n \in \omega: U_n \cap V\neq \emptyset\} \in p$ for every neighbourhood $V$ of $x$.
 
M. Hru\v{s}\'{a}k, I. Martínes-Ruiz and F. Hernandez-Hernandez in \cite{hrusak2007pseudocompactness} and V. Rodrigues, A. Tomita and Y. Ortiz-Castillo provided examples of spaces for which $X^\omega$ is pseudocompact and $\exp(X)$ is not. 
 
The question of whether there is a characterization for the spaces $X$ which have pseudocompact hyperspaces remains open. A natural class of spaces to explore this problem are the Isbell-Mrówka spaces, which are topological spaces associated to almost disjoint families. 
An almost disjoint family on a countable infinite set $N$ is an infinite set $\mathcal A$ of infinite subsets of $N$ which are pairwise almost disjoint, that is, for every $a, b\in \mathcal A$, $a\cap b$ is finite. A maximal almost disjoint family (on $N$), also called MAD family, is an almost disjoint family (on $N$) which is not contained in any other almost disjoint family. MAD families exist by Zorn's lemma and are uncountable. The least size of a MAD family is denoted by $\mathfrak a$. By $\mathcal I(\mathcal A)$ we denote the free ideal generated by $\mathcal A$, that is, $\mathcal I(\mathcal A)=\{X\subseteq \omega:\exists \mathcal B\in [\mathcal A]^{<\omega}\,X\subseteq^*\bigcup \mathcal B\}$, and $\mathcal I^+(\mathcal A)=\mathcal P(\omega)
  \setminus \mathcal I(\mathcal A)$

Given an almost disjoint family $\mathcal A$ on $\omega$, $\Psi(\mathcal A)$ is the set $\omega\cup\mathcal A$ topologized by the finest topology which makes $\omega$ open and discrete and each $a\subseteq \Psi(\mathcal A)$ is a sequence converging to the point $a\in \Psi(\mathcal A)$. A basis for this topology is 
$$\{\{n\}:n\in \omega\}\cup\{\{a\}\cup(a\setminus F):a\in \mathcal A\andd F\in [\omega]^{<\omega}\}.$$
The space $\Psi(\cA)$ is known as the Isbell-Mr\'owka space or $\Psi$-space of $\cA$. These spaces are Tychonoff, locally compact, zero-dimensional, separable, first-countable and not countably compact. We refer to the survey \cite{Hrusak2014} for more information on Isbell-Mrowka spaces. The following result is what makes psi spaces interesting for examining Ginsburg's questions.
	
\begin{proposition}
Let $\mathcal A$ be an almost disjoint family. Then $\mathcal A$ is a MAD family if and only if $\Psi(\mathcal A)$ is pseudocompact if and only if $\Psi(\mathcal A)^\omega$ is pseudocompact.\qed
\end{proposition}

As a shorthand, we say that an almost disjoint family $\mathcal A$ is pseudocompact if $\exp(\Psi(\mathcal A))$ is pseudocompact.
In \cite{hrusak2007pseudocompactness}, it is proved that $\mathfrak{p}=\mathfrak{c}$ implies that every MAD family is pseudocompact and that, on the other hand, $\mathfrak{h}<\mathfrak{c}$ implies the existence of a non-pseudocompact MAD family. It is also mentioned in \cite{hrusak2007pseudocompactness}, that whether $\exp(\Psi(\cA))$ is pseudocompact for every/some MAD family was explicitly asked by J. Cao and T. Nogura.
 
It is known that if $\exp(X)$ is pseudocompact, so is $X$ \cite{ginsburg1975some}. Thus, when restricting to Isbell-Mr\'owka spaces, the questions regarding the relations between the pseudocompactness of $X$ and of $\exp(X)$ boil down to the following:
	
	\begin{question}\label{question1}Is there a pseudocompact MAD family?
	\end{question}
	\begin{question}\label{question2}Is every MAD family pseudocompact?
	\end{question}
	
Question \ref{question2} had been partially answered in \cite{hrusak2007pseudocompactness} and was later solved in \cite{NewPaperPseudocompactness}, where it was proved that every MAD family is pseudocompact if and only if the Baire number of $\omega^*$ is strictly greater than $\mathfrak c$. In the last section of this paper we show that these statements are equivalent to every MAD family being $p$-pseudocompact for some free ultrafilter $p$.

Question \ref{question1} appears in \cite{Hrusak2014}, \cite{hruvsak2018pseudocompact} and \cite{hrusak2007pseudocompactness} and is still open.
	
	In \cite{intersecting} we introduced a new class of almost disjoint families called ``fin-intersecting''. 
Recall that a family $\cF\subseteq\cP(\omega)$ is centered if $\bigcap\mathcal B$ is infinite for every $\mathcal B\in[\mathcal F]^{<\omega}$.

\begin{definition}
A \emph{fin sequence} is a function $C:\omega\rightarrow [\omega]^{<\omega}\setminus \{\emptyset\}$ such that for all $n, m \in \omega$, if $n\neq m$ then $C(n)\cap C(m)=\emptyset$.

We say that an almost disjoint family $\mathcal A$ is \emph{fin-intersecting} if for every fin sequence $C$ there exists an infinite set $I\subseteq \omega$ such that $\{\{n \in I: a \cap C_n \neq \emptyset\}: a \in \mathcal A\}\setminus [I]^{<\omega}$ is centered.
	\end{definition}

In \cite{intersecting}, it was proved that every fin-intersecting MAD family is pseudocompact and that fin-intersecting MAD families exist under any of the following hypotheses:
\begin{itemize}
    \item $\mathfrak{a}<\mathfrak{s}$,
    \item $\mathfrak{ap}=\mathfrak{c}=\mathfrak{s}$
    \item we are in a model obtained by adding any amount of Cohen and/or Random reals over a model of \textsf{CH}.
\end{itemize}

In this paper, we provide new consistent examples for the existence of pseudocompact MAD families by proving that under $\mathfrak{ap}=\mathfrak c$ and $\diamondsuit(\mathfrak b)$ there exists such families, and discuss some properties regarding the possible structures for fin-intersecting MAD families by proving that under CH there exists a Cohen indestructible non fin-intersecting MAD family and by providing an example of a MAD family which contains a large (stationary-indexed) fin-intersecting subfamily but every club-indexed sub-family is not fin-intersecting.

Our notation is standard and we refer the reader to \cite{Blass2010} for cardinal invariants with the exception of $\mathfrak{ap}$ which may be found in \cite{Brendle} and will be defined below. For set theoretic notation, the reader may consult \cite{kunenNew} and for general topology we refer to \cite{engelking}.

\section{More models with pseudocompact MAD families.}

We recall that $\mathfrak{ap}$ is the least size of an almost disjoint family $\mathcal A$ that contains two subsets that cannot be weakly separated, that is, that contains a subfamily $\mathcal B$ such that there is no $X\subseteq \omega$ so that $|X\cap a|<\omega$ for every $a\in \mathcal A\setminus \mathcal B$ and $|X\cap a|=\omega$ for every $a\in \mathcal B$. In ZFC we have $\mathfrak p\leq \mathfrak{ap}\leq \mathfrak b$, and both inequalities are strictly consistent. For more on $\mathfrak{ap}$ the reader may consult \cite{Brendle}.

To simplify the notation, let $\pi_C(A)=\{n\in\omega:A\cap C(n)\neq\emptyset\}$, for a fin-sequence $C$ and an infinite set $A$. Let also $\pi_C[\cZ]=\{\pi_C(Z):Z\in\cZ\}$ and $E_C(X)=\{n\in\omega:C(n)\subseteq X\}$.

In the following, we use the fact that every family of less than $\mathfrak c$ infinite sets of natural numbers may be refined to an almost disjoint family \cite[Theorem 2.1]{mateweak}.

\begin{theorem}[$\mathfrak{ap=c}$]\label{theoremap}
There is a pseudocompact MAD family.
\end{theorem}

\begin{proof}
Denote the set of all fin-sequences by $\mathcal F$ and enumerate $\mathcal F\times [\omega]^\omega$ as $\{(C_\alpha, X_\alpha):\alpha<\mathfrak{c}\}$ in such a way that each pair appears cofinally. We may assume that $\{\bigcup_{n\in\omega}C_\alpha(n): \alpha<\omega\}$ is almost disjoint (this will guarantee that we can get to the step $\omega$ of the construction easily). 

We recursively define a family $\{a_\alpha: \alpha<\mathfrak c\}$ and families $\{\mathcal B_C^\alpha: \alpha<\mathfrak c\ \land\ C\in\mathcal{F}\}$ and $\{I_C^\alpha: \alpha<\mathfrak c \land\ C\in\mathcal{F}\}$ such that for every $\alpha<\mathfrak c$ and $C\in\cF$:

\begin{enumerate}
    \item\label{apitem1} $a_\alpha \in [\omega]^\omega$,
    
    \item $a_\beta\cap a_\alpha=^*\emptyset$ for $\beta<\alpha$,
    
    \item $\mathcal B_C^\alpha\subseteq \{a_\beta:\beta<\alpha\}$,

    \item\label{apitem4} $I_C^\alpha=\{\beta<\alpha: a_\beta \in \mathcal B_C^\alpha\}$,
    
    \item\label{apitem5} $\mathcal B_C^\beta\subseteq \mathcal B_C^\alpha$ for every $\beta<\alpha$,
    
    \item\label{apitem6} $\mathcal B_C^0=\emptyset$,
    
    \item\label{apitem7} $\mathcal B_C^\alpha=\bigcup_{\beta<\alpha}B_C^\beta$ if $\alpha$ is limit,
    
    \item $\{\pi_C(a_\beta): \beta \in I_C^\alpha\}$ is centered.
    
    \item\label{apitem9} If $\{\pi_C(a_\beta): \beta \in I_C^\alpha\}\cup \{\pi_C(a_\alpha)\}$ is centered then $a_\alpha \in \mathcal B_C^{\alpha+1}$, otherwise $\mathcal B^{\alpha+1}_C=\mathcal B^\alpha_C$.
    
    \item\label{apitem10} If:
    
    \begin{itemize}
    \item for every $a \in \mathcal B_{C_\alpha}^\alpha$ we have $a\subseteq^* X_\alpha$,
    \item $\{\pi_{C_\alpha}(a_\beta): \beta \in I_{C_\alpha}^\alpha\}\cup \{E_{C_\alpha}(X_\alpha)\}$ is not centered, and
    \item For every $k\geq 0$ and $\beta_0, \dots, \beta_{k-1}$, $\bigcup_{n \in J}C_\alpha(n)\setminus X_\alpha \in \mathcal I^+(\{a_\beta: \beta<\alpha\})$, where $J=\bigcap_{i<k}\pi_{C_\alpha}(a_{\beta_i})$.
    \end{itemize}
    
    Then $\{\pi_{C_\alpha}(a_\beta): \beta \in I_{C_\alpha}^\alpha\}\cup \{\pi_{C_\alpha}(a_\alpha)\}$ is centered and $a_\alpha \not \subseteq^* X_\alpha$.
\end{enumerate}

Assume that such construction is possible.
Let $\mathcal A=\{a_\alpha: \alpha<\mathfrak c\}$. We show $\mathcal A$ is pseudocompact. Fix $C\in \mathcal F$ and let $\mathcal B=\bigcup_{\alpha<\mathfrak c} \mathcal B_C^\alpha$ and $J_C=\{\beta<\mathfrak c:a_\beta\in \mathcal B\}$. We show that $\mathcal B$ is an accumulation point for $C$.

Let $U_0, \dots, U_l$ be open subsets of $\Psi(\mathcal A)$ such that $\mathcal B\in \langle U_0, \dots, U_l\rangle$. Let $X=\omega\cap \bigcup_{i\leq l} U_i$ and for each $i \leq l$, fix $\beta_i\in J_C$ and $m_i \in \omega$ such that $a_{\beta_i}\in \mathcal B\cap U_i$ and $a_{\beta_i}\setminus m_i\subseteq U_i$. Notice that for every $a \in \mathcal B$, $a\subseteq^* X$.
Let $\alpha>\beta_0, \dots, \beta_l$ be such that $(C, X)=(C_\alpha, X_\alpha)$.
For each $F\subseteq \mathfrak c$ and $\gamma<\mathfrak c$ let $G^F_\gamma=\bigcap\{\pi_{C_\gamma}(a_\beta):\beta\in F\}$. If $C_\gamma=C$ we simply write $G^F$. Let $F_0=\{\beta_i:i\leq l\}$.

     \textbf{Case 1:} $A=G^{F_0}\cap E_C(X)$  is infinite: in this case $C(n)\in \langle U_0, \dots, U_l\rangle$ for every $n \in A\setminus (\max\{m_0, \dots, m_l\})+1$, as needed.

        \textbf{Case 2:} $A=G^{F_0}\cap E_C(X)$ is finite.

        \textbf{Subcase 1:} for every finite $F\subseteq \alpha$, $\bigcup_{n \in G^F}C(n)\setminus X \in \mathcal I^+(\{a_\beta:\beta<\alpha\})$. In this case, $a_\alpha \in \mathcal B_C$ and $a_\alpha \not \subseteq^* X$ by (\ref{apitem10}), a contradiction.

        \textbf{Subcase 2:} There exists $F\subseteq \alpha$ such that $\bigcup_{n \in G^F}C(n)\setminus X \in \mathcal I(\{a_\beta:\beta<\alpha\})$. In this case, let $F_1$ be a finite subset of $\alpha$ for which there exists $K\subseteq \alpha$ of the smallest possible cardinality such that $\bigcup_{n \in G^{F_1}}C(n)\setminus X\subseteq^* \bigcup_{\beta \in K} a_\beta$. Fix an arbitrary $\beta_0 \in K$.

        \textbf{Claim:} $\beta_0 \in \mathcal B_\mathcal C^{\beta_0+1}$: we have to see that $\{\pi_C(a_\beta): \beta\in I_C^{\beta_0}\}\cup \{\pi_C(a_{\beta_0})\}$ is centered. Fix any finite $F_2\subseteq I_C^{\beta_0}$. Let $F=F_2\cup F_1\cup F_0$. Since $F_0\subseteq F$, we have that $G^F\cap E_C(X)$ is finite, thus $\bigcup_{n\in I_C^{F}}C(n)\setminus X$ is infinite. Since $F\supseteq F_1$, we have that $\bigcup_{n \in G^F}C(n)\setminus X\subseteq^* \bigcup_{\beta \in K} a_\beta$, thus $a_{\beta_0}\cap \bigcup_{n \in G^F}C(n) \subseteq a_{\beta_0} \cap\bigcup_{n \in G^{F_2}}C(n)$ is infinite by the minimality of $|K|$. This completes the proof of the claim.\\
        Thus $a_{\beta_0}\in\mathcal B$ but $a_{\beta_0}\cap(\omega\setminus X)$ is infinite, contradicting that $a\subseteq^* X$ for every $a\in\mathcal B$ which shows that Subcase 2 is not possible.

        \vspace{2mm}

    Now we verify that such construction is possible.    
    If $\alpha<\omega$, let $a_\alpha=\bigcup_{n\in \omega}C_\alpha(n)$. $B^\alpha_{C_\alpha}$ is defined as in (\ref{apitem6}) or (\ref{apitem9}) and $I^\alpha_C$ as in (\ref{apitem4}). It is easy, then, to check that (\ref{apitem1})-(\ref{apitem9}) hold and (\ref{apitem10}) will hold since $\pi_{C_\alpha}(a_\alpha)=\omega$, and since if $\{n\in \omega: C_\alpha\subseteq X_\alpha\}$ is not cofinite, then $a_\alpha=\bigcup_{n\in \omega} C_\alpha\setminus X_\alpha$ is infinite.

If $\omega\leq \alpha<\mathfrak c$ we must show how to construct $a_\alpha$ in the case where the hypothesis of (\ref{apitem10}) hold, otherwise $a_\alpha$ can be defined arbitrarily with the property that $|a_\alpha\cap a_\beta|<\omega$ for every $\beta<\alpha$ as $\alpha<\mathfrak{a}=\mathfrak{c}$. 
It also follows from $\mathfrak a=\mathfrak c$, that for each $F\subseteq \alpha$, there exists $a_F\subseteq \bigcup_{n\in G_{C_\alpha}^F}C_\alpha(n)\setminus X_\alpha$ almost disjoint with every $a_\beta$ for $\beta<\alpha$. 
The family $\{a_F: F\in [\alpha]^{<\omega}\}$ can be refined to an almost disjoint family $\{b_F: F\in [\alpha]^{<\omega}\}$ by the disjoint refinement lemma. Now use that $\mathfrak{ap}=\mathfrak c$ and let $a_\alpha$ be such that $a_\alpha \cap b_F$ is infinite for every $F \in [\alpha]^{<\omega}$ and $a_\alpha \cap a_\beta$ is finite for every $\beta<\alpha$. This completes the construction.
\end{proof}

Now we show that under $\diamondsuit(\mathfrak{b})$ there are pseudocompact MAD families of size $\omega_1$. Parametrized diamonds were introduced in \cite{parametrized} by M. Dzamonja, M. Hru\v{s}\'ak and J. Moore where many consequences of them are presented.
These principles are weakenings of Jensen's diamond principle $\diamondsuit$ but unlike it, parametrized diamond principles are consistent with the negation of \textsf{CH}. Recall that $\mathfrak{b}$ is the minimum size of a family $\mathcal B\subseteq\mathcal P(\omega^\omega)$ such that for every $f\in\omega^\omega$ there exists $b\in\mathcal B$ such that $b\nleq^* f$. Its parametrized diamond principle is then defined as follows:

$$\diamondsuit(\mathfrak{b})\equiv\ \forall F:2^{<\omega_1}\to \omega^\omega\ \textnormal{Borel}\ \exists g:\omega_1\to \omega^\omega\ \forall f\in2^{\omega_1}$$
$$\ \ \ \ \ \ \ \ \ \ \ \ \{\alpha\in\omega_1:F(f\rest\alpha) \ngeq^* g(\alpha)\}\textnormal{ is stationary},$$
where $F:2^{<\w_1}\rightarrow \omega^\omega$ is \emph{Borel}, if for every $\delta<\omega_1$ the restriction of $F$ to $2^\delta$ is a Borel map. 
For the next theorem, let us fix bijections $e_\alpha:\omega\to\alpha$ for every infinite ordinal $\alpha<\omega_1$.

\begin{theorem}
    If $\diamondsuit(\mathfrak{b})$ holds, there is a pseudocompact MAD family of size $\omega_1$.
\end{theorem}

\begin{proof}
    We can code the set $\omega_1\times([\omega]^\omega)^{\omega_1}\times2^{\omega_1}\times([\omega]^{<\omega})^\omega\times[\omega]^\omega)$ in $2^{\omega_1}$. More precisely, we may assume, by passing to a club subset of $\omega_1$ (e.g., the set of indecomposable ordinals), that every node $t\in2^{<\omega_1}$ codes a tuple 
    $$(\alpha,\{A_\beta:\beta<\alpha\},\sigma,C,X)$$
    such that:
    \begin{enumerate}
        \item $\alpha<\omega_1$,
        \item $\cA_\alpha=\{ A_\beta:\beta<\alpha\}\subseteq[\omega]^\omega$ is an almost disjoint family,
        \item $\sigma\in2^\alpha$,
        \item $C:\omega\to[\omega]^{<\omega}$ is a fin-sequence and
        \item $X\in[\omega]^\omega$.
    \end{enumerate}
\vspace{3mm}
We will use $\sigma\in2^\alpha$ to determine the subset $\cA_\alpha^\sigma=\{A_\beta\in\cA_\alpha:\sigma(\beta)=1\}$ of $\cA_\alpha$. 
Every branch $(\omega_1,\{A_\alpha:\beta<\omega_1\},f,C,X)$ in $2^{\omega_1}$ codes then an almost disjoint family $\{A_\alpha:\alpha<\omega_1\}$ along with a fin-sequence $C$, a subset $B\subseteq\{A_\alpha:\alpha<\omega_1\}$ given by the function $f\in2^{\omega_1}$ and an infinite set $X\in[\omega]^\omega$. 
The idea is to recursively construct an almost disjoint family $\cA=\{a_\alpha:\alpha<\omega_1\}$ along with a family $(B_C:C\textnormal{ is a fin-sequence})$ of elements of $\cP(\cA)$ such that $B_C$ is an accumulation point for every fin-sequence $C$. 
By considering $B_C$ as a characteristic function $f_C\in2^{\omega_1}$, we have that every such $5$-uple $(\alpha,\{a_\beta:\beta<\alpha\},f_C\rest\alpha,C,X)$ is coded by some node $t\in2^{<\omega_1}$. 
As mentioned before, by passing to the club of indecomposable ordinals, we may actually assume that $ht(t)=\alpha$. 
Hence, the guessing property of the $\diamondsuit(\mathfrak{b})$-sequence will help us to find a right choice of $a_\alpha$ at that step.\\

We define a function $F:2^{<\omega_1}\to\omega^\omega$ such that given $t=(\alpha,\{A_\beta:\beta<\alpha\},\sigma,C,X)$, we let $F(t)$ the constant function with value $0$ unless:
\begin{enumerate}[label=(\alph*)]
    \item\label{relevantitema} $\pi_C[\cA_\alpha^\sigma]$ is centered,
    \item\label{relevantitemb} $A\subseteq^* X$ for all $A\in\cA_\alpha^\sigma$ and
    \item\label{relevantitemc} $\pi_C[\cA_\alpha^\sigma]\cup\{E_C(X)\}$ is not centered.
\end{enumerate}
 We say that $t\in2^{<\omega_1}$ is \emph{relevant} if it satisfies the previous requirements. 
 Notice that these conditions are Borel. For example, item (a) can be written as:
 $$\forall G\in[\alpha]^{<\omega}\ \big(\forall \gamma\in G\  (\sigma(A)=1)\Rightarrow\exists^\infty n\in\omega\ \forall A\in G(A\cap C(n)\neq\emptyset)\big),$$
and similar statements describe (b) and (c).
Fix a relevant $t\in2^\alpha$ coding a tuple $(\alpha,\{A_\beta:\beta<\alpha\},\sigma,C,X)$ and define 
$$G_n^t=\{A_\beta\in\cA_\alpha:\exists m\leq n\ (e_\alpha(m)=\beta \land \sigma(\beta)=1)\}.$$
Also, for every $n\in\omega$ let $I_n=\bigcap_{A\in G_n^t}\pi_C(A)$.
Thus there is $N(t)\in\omega$ such that $I_{N(t)}\cap E(X)$ is finite by item \ref{relevantitemc}, which implies that $C(m)\setminus X\neq\emptyset$ for all but finitely many $m\in I_{N(t)}$. 
Henceforth $Y_k^t=\bigcup_{m\in I_k} C(m)\setminus X$ is infinite for $k=N(t)$. Since the family $\{G_n^t:n\in\omega\}$ is increasing, we get by item \ref{relevantitema} that $\{I_n:n\in\omega\}$ is a decreasing family of infinite sets and then $\{Y_k^t:k\in\omega\}$ is also a decreasing family consisting of infinite sets.

If there is a $k\in\omega$ such that $Y_k^t\in\cI(\cA_\alpha)$, define $F(t)(k)=0$ for every $k\in\omega$. Otherwise, we have that
\begin{itemize}
    \item[\ding{72}]\label{blackstar} $\forall k\in\omega\ \big(Y_k^t\in\cI^+(\cA_\alpha)\big)$,
\end{itemize}
which is also a Borel condition.
In the later case define 
$$F(t)(k)=\min\bigg(Y_k^t\big\backslash\Big(\bigcup_{i<k}A_{e_\alpha(i)}\Big)\bigg).$$\\

Let $g\in\omega^\omega$ be a $\diamondsuit(\mathfrak{b})$-sequence for $F$. We recursively define $\cA=\{a_\alpha:\alpha<\omega_1\}$ by
$$a_\alpha=\bigcup_{n\in\omega}\bigg(g_\alpha(n)\big\backslash\Big(\bigcup_{i<n}a_{e_\alpha(i)}\Big)\bigg)$$

where $g_\alpha$ stands for $g(\alpha)$. Clearly $\{a_\alpha:\alpha<\omega_1\}$ is almost disjoint and we can assume, by increasing each $g_\alpha$ if necessary, that each $a_\alpha$ is infinite. Given a fin-sequence $C$, we define $B_C\subseteq\cA$ recursively by defining $B_C^\emptyset=\emptyset$, $B_C^\alpha=\bigcup_{\beta<\alpha}B_C^\beta$ if $\alpha$ is limit, and $B_C^\alpha=B_C^\beta\cup\{a_\alpha\}$ if $\alpha=\beta+1$ and $\{\pi(a):a\in B_C^\beta\}\cup\{a_\alpha\}$ is centered. If the later set is not centered, let $B_C^\alpha=B_C^\beta$. Thus $f_C\in2^{\omega_1}$ given by $f_C(\alpha)=1$ if and only ff $a_\alpha\in B_C$, codes $B_C$, the candidate for an accumulation point for $C$.\\

We claim that $\cA$ is pseudocompact. 
Given a fin-sequence $C$ and an open neighborhood $\langle U_0,\ldots, U_n\rangle$ of $B_C$, let $X=\bigcup_{i\leq n}U_i$. We want to show that there are infinitely many $n\in\omega$ such that $C(n)\subseteq X$ and $C(n)\cap U_i\neq\emptyset$ for every $i\leq n$.
Note that $a\subseteq^* X$ for every $a\in B_C$ and that
for every $i\leq n$, we can find $a_{\alpha_i}\in U_i\cap B_C$ (equivalently $\{a_{\alpha_i}\}\in U_i^-$), so $a_{\alpha_i}\subseteq^*U_i$. 
As shrinking $U_i$ makes our work harder, we may assume, without loss of generality, that for every $i\leq n$ there is $l_i\in\omega$ such that $U_i=\{a_{\alpha_i}\}\cup(a_{\alpha_i}\setminus l_i)$.\\
By definition, $\{\pi(a_\alpha):a_\alpha\in B_C\}$ is centered, in particular $\big|\bigcap_{i\leq n}\pi_C(a_{\alpha_i})\big|=\omega$. 
Take $l\in\omega$ such that $C(m)\cap\bigcup_{i\leq n}(l_i+1)=\emptyset$ for every $m\geq l$.
Let $J=\Big(\bigcap_{i\leq n}\pi_C(a_{\alpha_i})\Big)\backslash l$, and notice that for any $m\in J$ we have that $C(m)\in U_i^-$ for every $i\leq n$. Thus if $J\cap E_C(X)$ is infinite we are done, as $C(m)\in\langle U_0,\ldots,U_n\rangle$ for every $m\in J\cap E_C(X)$. It remains to show that it is indeed the case.

Assume otherwise that $J\cap E_C(X)$ is finite. Then we have that 
\begin{itemize}
    \item[\ding{91}]\label{pencil} $\big|\bigcap_{i\leq n}\pi_C(a_{\alpha_i})\cap E_C(X)\big|<\omega$.
\end{itemize} 
Consider $b=(\omega_1,\{a_\alpha:\alpha<\omega_1\},f_C,C,X)$ being a branch in $2^{\omega_1}$ given by our coding. We can find an $\alpha>\max\{\alpha_i:i\leq n\}$ such that the $\diamondsuit(\mathfrak{b})$-sequence guesses $b$ at $\alpha$, that is, $g_\alpha\nleq^* F(b\rest\alpha)$, with
$$b\rest\alpha=(\alpha,\{a_\beta:\beta<\alpha\},f_C\rest\alpha,C,X).$$
It is clear that $b\rest\alpha$ satisfies properties (1)-(5). With $\cA_\alpha=\{a_\beta:\beta<\alpha\}$, $\sigma=f_C\rest\alpha$ and $\cA_\alpha^\sigma=\{a_\beta:\beta<\alpha\land\sigma(\alpha=1)\}$ it also follows from the definition of $B_C$ that $b\rest\alpha$ satisfies \ref{relevantitema} and \ref{relevantitemb}. Moreover, from \ding{91}, we get that $b\rest\alpha$ satisfies \ref{relevantitemc}.

\vspace{5mm}
\underline{\textit{Claim:}} Fix $t=b\rest\alpha$, then $Y_k^t\in\cI^+(\cA_\alpha)$ for every $k\in\omega$.\\
\underline{\textit{Proof of Claim:}} Assume otherwise that there is $k_0\in\omega$ such that $Y_{k_0}^t\in\cI(\cA_\alpha)$. As $\{Y_k^t:k\in\omega\}$ is a decreasing sequence, that means that $Y_k^t\in\cI(\cA_\alpha)$ for every $k\geq k_0$. 
Let $\{a_{\beta_0},\ldots,a_{\beta_r}\}\in[\cA_\alpha]^{<\omega}$ be such that $Y_{k_0}^t\subseteq^*\bigcup_{j\leq r}a_{\beta_j}$ and each $a_{\beta_j}$ has infinite intersection with $Y_{k_0}^t$. As $Y_{k_0}^t$ is disjoint from $X$ and every $a\in B_C$ is almost contained in $X$, we have that $\sigma(a_{\beta_j})=0$ for every $j\leq r$, i.e., $a_{\beta_j}\notin B_C$. 
As every $\beta_j$ is smaller than $\alpha$, this only happens if there is a finite set $G_j\in[\cA_\alpha^\sigma]^{<\omega}$ such that $|\pi_C(a_{\beta_j})\cap\bigcap_{a\in G_j}\pi_c(a)|<\omega$. 
Let $N>k_0$ and such that $G_j\subseteq G_N^t$ for every $j\leq N$, then clearly $|\pi_C(a_{\beta_j})\cap I_N^t|<\omega$ for every $j\leq r$ (recall that $I_N^t=\bigcap_{a\in G_N^t}\pi_c(a)$). 
This implies that $a_{\beta_j}\cap\bigcup_{m\in I_n^t}C(m)$ is finite for every $j\leq r$.
On the other hand
$$\bigcup_{m\in I_N^t}C(m)\subseteq\bigcup_{m\in I_N^t}C(m)\setminus X=Y_N^t\subseteq Y_{k_0}^t\subseteq^*\bigcup_{j\leq r}a_{\beta_j}$$
which is a contradiction.\hfill$_{Claim}\square$

\vspace{5mm}
By the claim, we have that $t$ satisfies \ding{72} and then $F(t)(k)=\min(Y_k^t\setminus(\bigcup_{i<k}a_{e_\alpha(i)}))$. As $g_\alpha\nleq^* F(t)$, there are infinitely many $k\in\omega$ such that $g_\alpha(k)>F(t)(k)$. 
In particular, $a_\alpha\cap Y_k^t$ is infinite for every $k\in\omega$ since $\{Y_k^t:k\in\omega\}$ is decreasing. 
Finally, given $H\in[B_C\rest\alpha]^{<\omega}$, there is $s\in\omega$ such that $H\subseteq G_s^t$ and in consequence $\bigcap_{a\in H}\pi_C(a)\subseteq I_s^t$. 
Given that $a_\alpha$ has infinite intersection with $Y_s^t=\bigcup_{m\in I_s^t}C(m)\setminus X$, from which it is easily seen that $I_s^t\cap\pi_C(a_\alpha)$ is infinite. 
Therefore $\{\pi_C(a):a\in\cA_\alpha\}\cup\{\pi_C(a_\alpha)\}$ is centered and $a_\alpha\in B_C$, but $a\nsubseteq^* X$ because each $Y_k^t$ is disjoint from $X$, contradicting that $a\subseteq^* X$ for every $a\in B_C$.
Thus, assuming \ding{91} leads us to a contradiction, which implies that $J\cap E_C(X)$ is infinite.
\end{proof}

The previous result and Theorem \ref{theoremap} follow the same idea of recursively constructing the almost disjoint family and approximations for all accumulation points needed to show that the hyperspace is pseudocompact. A different idea was used in \cite{intersecting} to show that if $\mathfrak{ap}=\mathfrak c$ and also $\mathfrak s=\mathfrak c$ then there is even a fin-intersecting MAD family.
We do not know if the two results presented here can be strengthened in order to prove the existence of fin-intersecting MAD families:

\begin{question}
    Does it follow from either $\mathfrak{ap}=\mathfrak{c}$ or $\diamondsuit(\mathfrak{b})$ that there are fin-intersecting MAD families?
\end{question}

\section{On the structure of non fin-intersecting MAD families}

Fin-intersecting almost disjoint families were introduced in \cite{intersecting}. Every fin-intersecting MAD family is pseudocompact, and this is a stronger property than pseudocompactness: the existence of a non pseudocompact MAD family is independent of ZFC, yet non fin-intersecting MAD families do exist in ZFC. It is not known if non pseudocompact MAD families of size $\omega_1$ exist, but non fin-intersecting MAD families of size $\omega_1$ exist under $\mathfrak{s}=\omega_1$.

In this section, we discuss how badly the fin-intersecting property may fail.
First, we show that under CH, there exists a non fin-intersecting MAD family whose non fin-intersectingness is preserved by adding Cohen reals.

 \begin{theorem}[CH] There exists a non fin-intersecting MAD family $\mathcal A$ such that, $\mathbb C \Vdash ``\check{\mathcal A} \text { is MAD and non fin-intersecting}$''. More specifically, there exists such a MAD family for which there exists a fin sequence $C$ which witnesses the non fin-intersectness of $\mathcal A$ in every (iterated, with finite supports) Cohen extension.
    \end{theorem}

 \begin{proof} 
Enumerate $\mathbb C\times \{\tau \in V^{\mathbb C}: \tau \text{ is a nice name for a subset of } \check \omega\}$ as $\{(p_\alpha, \tau_\alpha): \alpha \in [\omega, \omega_1)\}$.
 Let $\mathcal A$ be an infinite countable almost disjoint family. Enumerate it as $\{a_\alpha: \alpha<\omega\}$.  We will expand it to a MAD family as intended. Let $C:\omega\rightarrow [\omega]^{<\omega}\setminus \{\emptyset\}$ be a fin sequence such that for every $n\in \omega$ and $\alpha<\omega$, $|C(n)|\geq n+1$ and $|a_\alpha\cap C(n)|\leq 1$.

 We recursively define, for every $\alpha \in [\omega, \omega_1)$, sets $a_\alpha, b_\alpha, c_\alpha \in [\omega]^\omega$ and $\mathcal A_\alpha$ such that:

\begin{enumerate}[label=\alph*)]
\item $\mathcal A_\alpha=\{a_\xi: \xi < \alpha\}\cup\{b_\xi, c_\xi: \omega\leq \xi< \alpha\}$ is an almost disjoint family.
\item for every $x \in \mathcal A_\alpha$ and $n \in \omega$, $|x\cap C(n)|\leq 1$.
\item If $p_\alpha \Vdash ``|\tau_\alpha|=\omega$'' and $p_\alpha \Vdash ``\tau_\alpha \notin \cI(\check{\mathcal A_\alpha})$'' then $p_\alpha \Vdash  ``|\check a_\alpha \cap\tau_\alpha|=\omega$''.
\item If $p_\alpha \Vdash ``|\tau_\alpha|=\omega$'' then $p_\alpha \Vdash ``|\{n \in \tau_\alpha: \check C(n)\cap \check b_\alpha\}|=|\{n \in \tau_\alpha: \check C(n)\cap \check c_\alpha\}|=\omega$''.
\item If $p_\alpha \Vdash ``|\tau_\alpha|=\omega$'' then $p_\alpha \Vdash ``\{n \in \omega: \check C(n)\cap  b_\alpha\}\cap\{n \in \omega: \check C(n)\cap  c_\alpha\}=\emptyset$''.
\end{enumerate}

We show how to define $a_\alpha, b_\alpha, c_\alpha \in [\omega]^\omega$ having defined all these objects for all $\xi \in [\omega, \alpha)$. We define $a_\alpha$ and then $b_\alpha, c_\alpha$.

\textbf{Case 1:} Defining $x=a_\alpha$, $b_\alpha$, $c_\alpha$ in case the hypothesis c) (for $a_\alpha$) or d) (for $b_\alpha$, $c_\alpha$) fails. In this case, $x$ only needs to satisfy $b)$ and be almost disjoint with the elements of the almost disjoint family that have been constructed so far. Enumerate the previously constructed elements as $\{d_n: n \in \omega\}$, choose $x_n \in C(n)\setminus \bigcup_{i<n}d_i$ and let $x=\{x_n: n \in \omega\}$. Note that this is possible by the assumption on the size of $C(n)$.

\textbf{Case 2:} Defining $a_\alpha$ in case $p_\alpha \Vdash ``|\tau_\alpha|=\omega$'' and $p_\alpha \Vdash ``\tau_\alpha \notin I(\check{\mathcal A_\alpha})$''.
Enumerate $\omega\times \{q \in \mathbb C: q\leq p_\alpha\}=\{(l_i, q_i): i\in\omega\}$. Enumerate also $\mathcal A_\alpha$ as $\{d_n: n \in \omega\}$. For each $i\in \omega$, let $r_i\leq q_i$ and $m_i \in \omega$ be such that $m_i\geq l_i$ and $r_i \Vdash ``\check m_i \in \tau_\alpha \setminus \check{\bigcup_{j<i}d_j}$''. 
Now let $a_\alpha=\{m_i: i \in \omega\}$. To see that $a_\alpha$ satisfies $c$, if it does not we would have $q\leq p_\alpha$ and $l\in \omega$ such that $q \Vdash``\check a_\alpha \cap \tau_\alpha \subseteq \check l$''. Let $i$ be such that $(l, q)=(l_i, q_i)$. Then $r_i\leq q$ and $r_i ``\Vdash \check m_i \in a_\alpha \cap \tau_\alpha \setminus \check l_i$'', a contradiction.

\textbf{Case 3:} Defining $b_\alpha$, $c_\alpha$ in case $p_\alpha \Vdash`` |\tau_\alpha|=\omega$''. 
Enumerate $\omega\times \{q \in \mathbb C: q\leq p_\alpha\}=\{(l_i, q_i): q\leq p_\alpha\}$. Enumerate $\mathcal A_\alpha\cup \{a_\alpha\}$ as $\{d_n: n \in \omega\}$. For each $i\in \omega$, let $r_i\leq q_i$, $m^b_i \in \omega$ and $m^c_i \in \omega$ be such that $m^b_i\geq l_i$, $m^c_i\geq l_i$ and $r_i \Vdash \check m^b_i, m^c_i \in \tau_\alpha$. We can choose these numbers so that the family $(m^b_i, m^c_i: i \in \omega)$ is injective. Now for each $i$, let $z_i^b \in C(m_i^b)\setminus \bigcup_{j<l_i} d_i$, $z_i^c \in C(m_i^c)\setminus \bigcup_{j<l_i} d_i$. Finally, let $b_\alpha=\{z^b_i: i \in \omega\}$ and $c_\alpha=\{z^c_i: i \in \omega\}$. It is clear that d) and e) hold.

\textbf{The construction works:} Let $\mathcal A'=\{a_\xi: \xi<\omega_1\}\cup\{b_\xi: \omega\leq \xi<\omega_1\}\cup\{c_\xi: \omega\leq \xi<\omega_1\}$. We first observe that $\mathcal A$ is MAD and non-fin-intersecting if we force with $\mathbb C$. By a), $\mathcal A$ is an almost disjoint family. By c), $\mathcal A$ is forced to be MAD, and by d) and e), $\mathcal A$ is forced to be non-fin-intersecting.

For the iteration, given a cardinal $\kappa$, we shall show that $C$ witnesses that $\mathbb C_\kappa \Vdash ``\check{\mathcal A} \text { is non-fin-intersecting}$'': Let $G$ be $\mathbb C_\kappa$-generic over $V$. Fix $I \in [\omega]^\omega \cap V[G]$. There exists a $\mathbb C$-generic filter $H$ over $V$ such that $I\in V[H]\subseteq V[G]$. Thus, in $V[H]$, $\{n \in I: C(n)\cap b_\alpha\}$ and $\{n \in I: C(n)\cap c_\alpha\}$ are infinite and disjoint, but these properties are absolute.
\end{proof}

We know that every countable (of size less than $<\mathfrak s$) almost disjoint family is fin-intersecting. Thus, it is natural to ask whether it is consistent that there exists a MAD family whose every uncountable subfamily fails to be fin-intersecting (e.g., assuming \textsf{CH}).

\begin{question}\label{question9}
    Is it consistent with \textsf{CH} that there is an almost disjoint family such that every uncountable subfamily fails to be fin-intersecting?
\end{question}

A tight MAD family, also called $\aleph_0$-MAD family, is an almost disjoint family $\mathcal A$ such that for every sequence $(X_n: n\in \mathbb N)$ of elements of $\mathcal I^+(\mathcal A)$ there exists $a \in \mathcal A$ such that $a\cap X_n$ is infinite for every $n \in \omega$ (notice that equivalently we may only ask for $a\cap X_n\neq\emptyset$ for every $n\in\omega$). Tight MAD families are Cohen-indestructible and every Cohen-indestructible MAD family is ``somewhere tight'' \cite{kurilic2001cohen}. Thus, the following result gives us an example of a Cohen indestructible MAD family in which certain (very large) subfamilies are not fin-intersecting. 
This complements with the previous result, which gave us a Cohen-indestructible MAD family whose non-fin-intersectingness was preserved.

\begin{theorem}[CH] There exists a tight MAD family $\{a_\alpha: \alpha<\omega_1\}$ such that for every club $K\subseteq \omega_1$, $\{a_\alpha: \alpha\in K\}$ is not fin-intersecting and there exists a stationary set $S$ such that $\{a_\alpha: \alpha\in S\}$ is fin-intersecting.
\end{theorem}

\begin{proof} Enumerate $[\omega]^{\omega}=\{I_\beta: 1\leq \beta<\omega_1\}$ and partition $\omega_1$ into $\omega_1$ stationary sets, $\{S_\alpha: \alpha<\omega_1\}$. We assume $0 \notin S_0$. Given $\beta$ in $S_0$, let $\beta^+$ be the successor of $\beta$ in $S_0$.

Let $C$ be a fin-sequence such that for every $n \in \omega$, $|C(n)|\geq n+1$ and that $\bigcup_{n\in \omega}C(n)=\omega$.
Let $\{D_\alpha: \alpha\in S_0\}$ enumerate all fin-sequences and enumerate $([\omega]^\omega)^\omega=\{X_\alpha: \alpha \in S_0\}$.

We recursively construct an almost disjoint family $\{a_\alpha: \alpha<\omega_1\}$ and for each $\beta\in S_0$, families $(I_\beta^\alpha: \alpha \in S_0\setminus \beta)$ satisfying that for every $\alpha$:

\begin{enumerate}[label=\alph*)]
\item For every $n \in \omega$, $|a_\alpha\cap C(n)|\leq 1$.
\item $a_\alpha\cap a_\beta$ is finite for every $\beta<\alpha$.
\item If $\alpha \in S_\beta$ and $\beta>0$, $a_\alpha\subseteq \bigcup_{n \in I_\beta} C_n$.
\item If $\alpha \in S_0$ and for every $n \in \omega$, $X_\alpha(n)\in \mathcal I^+(\{a_\xi: \xi<\alpha\})$ then for every $n \in \omega$, $X_\alpha(n)\cap a_\alpha\neq \emptyset$.
\item The family $\{I_\beta^\alpha: \alpha \in S_0\setminus \beta\}$ is $\subseteq^*$-decreasing for every $\beta\in S_0$, and each element of this sequence is infinite or empty.
\item If $\alpha \in S_0$, let $I_\alpha^\alpha$ be an infinite set such that for every $\xi<\alpha$, $\{n \in I_\alpha^\alpha: C_\alpha(n)\cap a_\xi\}$ is either finite or cofinite in $I_\alpha^\alpha$.
\item If $\alpha, \beta \in S_0$, $\beta\leq \alpha$ and $\bigcup_{n \in I_\beta^\alpha}C_\beta(n)\in \mathcal I^+(\{a_\xi:\xi<\alpha\})$, then $I_\beta^{\alpha^+}=\{n\in I_\beta^{\alpha}:C_\beta(n)\cap a_{\alpha}\neq \emptyset\}$ is infinite.
\item If $\alpha, \beta \in S_0$, $\beta<\alpha$,  $\alpha$ is limit in the order type of $S_0$ and for every $\xi\in [\beta, \alpha)\cap S_0$, $I_\beta^\xi$ is infinite, then so is $I_\beta^\alpha$.

\end{enumerate}
 
At step $\alpha$, we break the construction into two cases:

\textbf{Case 1:} there exists $\beta>0$ such that $\alpha \in S_\beta$. If the hypothesis for $c)$ hold and if $\alpha\neq0$ fix a surjection $\phi: \omega\rightarrow \alpha$.In case $\alpha=0$ let $\phi:\omega\to\{\emptyset\}$ (that is, $\phi$ is the constant function with value $\emptyset$). Write $I_\alpha=\{k_n: n \in \omega\}$ in increasing order and  $x_n \in C_{k_n} \setminus \bigcup_{i<n} a_{\phi(i)}$. Then let $a_\alpha=\{x_n: n\in \omega\}$.

\textbf{Case 2:} $\alpha \in S_0$ and then $\alpha>0$. Let $\phi: \omega\rightarrow \alpha$ be onto. Partition $\omega$ into infinite sets, $\omega=\bigcup_{j \in \omega} P_j$. If the hypothesis for d) holds, recursively define $x_n$ and $k_n$ such that $k_n$ is strictly increasing and whenever $n \in P_j$, $x_n \in C_{k_{2n}}\cap (X_i\setminus \bigcup_{i<n}a_i)$, then let $a'_\alpha=\{x_n: n\in \omega\}$. If the hypothesis for d) does not hold, let $a'_\alpha$ be almost disjoint with every $a_\xi$ with $\xi<\alpha$.

If there is no $\beta\leq \alpha$ such that $\bigcup_{n\in I_\beta^\alpha}C_\beta(n)\in \mathcal I^+(\{a_\xi: \xi<\alpha\})$, then we let $a_\alpha=a'_\alpha$. 
Otherwise, enumerate these $\beta$'s as $\{\beta_n: n\in \omega\}$ and select a strictly increasing sequence $\{k_n':n\in\omega\}$ disjoint from $\{k_{2n}:n \in \omega\}$ and $x'_n$ such that whenever $n\in P_j$, $k'_n\in I_{\beta_j}^\alpha$ and, $x'_n \in C(k'_n)\setminus \bigcup_{i<n}a_{\phi(i)}$.
Let $a_\alpha''=\{x'_n: n \in \omega\}$ and finally define $a_\alpha=a_\alpha'\cup a_\alpha''$.

The definition of $I_\alpha^\alpha$follows easily from the fact that $\mathfrak s=\mathfrak c$ and in order to satisfy 
$h)$, we just need to take a pseudointersection $I_\beta^\alpha$.
This finishes the construction. Let $\mathcal A=\{a_\alpha: \alpha<\omega_1\}$.

\textbf{$\mathcal A$ is tight: } Let $X:\omega\rightarrow \mathcal I^+(\mathcal A)$. Let $\alpha \in S_0$ be such that $X=X_\alpha$. Then $a_\alpha\cap X(n)\neq\emptyset$ for every $n \in \omega$.

\textbf{$\mathcal A$ has no fin-intersecting subfamily indexed by a club } Let $K$ be a club. We show that $\{a_\alpha: \alpha\in K\}$ is not fin-intersecting. Fix $I\in [\omega]^\omega$. Partition it into two disjoint infinite subsets, which are enumerated as $I_\beta$ and $I_{\beta'}$ with $\beta, \beta'>0$. Let $\alpha \in K\cap S_\beta$ and $\alpha' \in K\cap S_{\beta'}$. Then $a_\alpha \subseteq \bigcup_{n \in I_\beta}C(n)$ and $a_{\alpha'} \subseteq \bigcup_{n \in I_{\beta'}}C(n)$, which implies that $\{n \in I: C(n)\cap a_\alpha\neq \emptyset \}$ and $\{n \in I: C(n)\cap a_{\alpha'}\neq \emptyset \}$ are infinite and disjoint. Thus, $C$ witnesses that $\{a_\alpha: \alpha\in K\}$ is not fin-intersecting.

\textbf{The fa $\{a_\xi: \xi \in S_0\}$ is fin-intersecting:} given a fin-sequence $D$, there exists $\beta \in S_0$ such that $D=D_\beta$. 
If there exists $J\in [\omega]^\omega$ such that $\bigcup_{n\in J}D_\beta(n)\in \mathcal I(\mathcal A)$, by letting $\mathcal A'\subseteq \mathcal A$ be finite and minimal such that there exists an infinite $J'\subseteq J$ such that $\bigcup_{n\in J'} C_\beta(n)\subseteq^* \bigcup \mathcal A'$, it can be easily observed that $\mathcal A'=\{a\in \mathcal A:|\{n \in J': a\cap C_\beta(n)\neq \emptyset\}|=\omega\}$ and that $J'=^*\bigcap_{a\in \mathcal A'}\{n \in J': a\cap C_\beta(n)\neq \emptyset\}$. Hence $J'$ witnesses that fin-intersectness for $D$.

Thus, we may assume that $\bigcup_{n\in J}C_\beta(n)\in \mathcal I^+(\mathcal A)$ for every $J\in [\omega]^\omega$. In particular $\bigcup_{n \in I_\beta^\alpha}C_\beta(n)\in \mathcal I^+(\mathcal A)$ for every $\alpha \in S_0\setminus \beta$ as every $I_\beta^\alpha$ is infinite for $g)$.

Let $F\subseteq S_0$ be finite. By f), $J_0=\{\{n \in I_\beta^\beta: C_\beta(n)\cap a_\alpha\neq \emptyset\}:\alpha \in F\cap \beta\}\setminus [\omega]^{<{\omega}}$ is centered. 
Write $F\setminus \beta=\{\alpha_0, \dots, \alpha_{l-1}\}$ in increasing order. If $l=0$ we are done. If not, by finite induction using g) we can show, for $i<l$, that $J_{i+1}=\{n\in J_i: a_{\alpha_i}\cap C_\beta(n)\neq \emptyset\}$ is infinite, which completes the proof.
\end{proof}

The following question first posed in \cite{intersecting} in a different way, remains open and can be naturally combined with Question \ref{question9}. 

\begin{question}
    Is $\mathfrak s$ characterized by being the first cardinal $\kappa$ for which there is a non fin-intersecting almost disjoint family of size $\kappa$?
\end{question}

\begin{question}
    If the answer for the previous question is positive: Is there a non-fin-intersecting MAD family such that each subfamily of size $\mathfrak s$ is not fin-intersecting?
\end{question}
\section{The \texorpdfstring{$p$}{p}-pseudocompactness of Isbell-Mrówka spaces and their hyperspaces}

If $X$ is a topological space with no isolated points, the Baire number of $X$, or Novák's number of $X$, denoted by $n(X)$, is the least cardinality of a collection of open dense subsets of $X$ with empty intersection. The cardinal $n(\omega^*)$ was defined and studied in \cite{balcar1980space}. In ZFC, $\mathfrak p^+, \mathfrak h\leq n(\omega^*)\leq 2^\mathfrak c$ and all inequalities are consistently strict.

In \cite{NewPaperPseudocompactness}, it was proved that $n(\omega^*)>\mathfrak c$ is equivalent to every MAD family being pseudocompact. In this section, we modify the arguments in that paper to show that these statements are equivalent to every MAD family being $p$-pseudocompact for some free ultrafilter $p$.

In \cite{salvador}, S. García-Ferreira defined the notions of $(\kappa, A)$ pseudocompactness, where $\kappa$ is a cardinal and $A\subseteq \omega^*$.  A space $X$ is said to be $(\kappa, A)$-pseudocompact if and only if for every family of $\kappa$ sequences $(U^\alpha: \alpha<\kappa)$ of open sets there exists $p \in \omega^*$ such that for every $\alpha<\kappa$, $U^\alpha$ has a $p$-limit point. If $X$ is Tychonoff, then for every $\alpha\leq \omega$, $(\alpha, \omega^*)$-pseudocompactness is equivalent to $X^\alpha$ being pseudocompact. For a proof and more on these notions, see \cite{salvador} or \cite{salvador2}. Clearly, $p$-pseudocompactness is equivalent to the space being $(\kappa, \{p\})$-pseudocompact for every $\kappa$.

\begin{proposition}\label{cappseudocompact}
Let $A\subseteq\omega^*$ and $\mathcal A$ be an almost disjoint family. The following are equivalent:

\begin{enumerate}
    \item\label{itema} $\exp(\Psi(\mathcal A))$ is $p$-pseudocompact for some $p \in A$,
    \item\label{itemb} $\Psi(\mathcal A)$ is $p$-pseudocompact for some $p \in A$,
    \item\label{itemc} $\Psi(\mathcal A)$ is $(\mathfrak c, A)$-pseudocompact for some $p\in A$.

\end{enumerate}
\end{proposition}

\begin{proof}
Firstly, recall that, since $\Psi(\mathcal A)$ is Tychonoff, if $p \in \omega^*$, by \cite[Theorem 2.4.]{ginsburg1975some}, the $p$-pseudocompactness of $\Psi(\mathcal A)$ is equivalent to the $p$-pseudocompactness of $\exp(\Psi(\mathcal A))$. Thus, (\ref{itema}) and (\ref{itemb}) are equivalent.

That (\ref{itemb}) implies (\ref{itemc}) is trivial. For (\ref{itemc}) implies (\ref{itemb}) we proceed as follows:

For every $f \in \omega^\omega$, let $U_f$ be the sequence of open subsets of $\Psi(\mathcal A)$ given by $U_f(n)=\{f(n)\}$. Let $p \in A$ be such that every $U_f$ has a $p$-limit point. Then $\Psi(\mathcal A)$ is $p$-pseudocompact. Given a sequence of nonempty open sets  $V=(V_n: n \in \omega)$, by the density of $\omega$ there exists $f \in \prod_{n \in \omega} V_n$. It is clear then that a $p$-limit point of $U_f$ is a $p$-limit point of $V$.
\end{proof}

By setting $A=\{p\}$ we obtain the following corollary:

\begin{corollary}\label{cappseudocompact2}
Let $p \in \omega^*$ and $\mathcal A$ be an almost disjoint family. The following are equivalent:

\begin{enumerate}
    \item\label{itema2} $\exp(\Psi(\mathcal A))$ is $p$-pseudocompact,
    \item\label{itemb2} $\Psi(\mathcal A)$ is $p$-pseudocompact,
    \item\label{itemc2} $\Psi(\mathcal A)$ is $(\mathfrak c, \{p\})$-pseudocompact.

\end{enumerate}
\end{corollary}

\begin{proposition}
For every $\kappa<n(\omega^*)$ and for every MAD family $\mathcal A$, there exists $p \in \omega^*$ such that $\Psi(\mathcal A)$ is $(\kappa, \omega^*)$-pseudocompact.
\end{proposition}
\begin{proof}
   Let $((U^\alpha(n):n\in\omega): \alpha<\kappa)$ be such that each $(U^\alpha(n):n\in\omega)$ is a sequence of nonempty open subsets of $\Psi(\mathcal A)$. For each $\alpha<\kappa$, let $f_\alpha \in \prod_{n<\omega}(U^\alpha(n)\cap\omega)$ and let 
    $$D_\alpha=\{p \in \omega^*: \exists B\in p\ (f_\alpha|B \text{ is constant }) \vee (f_\alpha|B \text{ is injective and } \exists a \in \mathcal A \, f_\alpha[B]\subseteq a)\}.$$

    For each $\alpha<\kappa$, $D_\alpha$ is open and dense:
    
    \begin{itemize}
        \item It is open: given $p \in D_\alpha$, let $B\in p$ be a witness for $p \in D_\alpha$. Then $B^*=\{q \in \omega^*: B\in q\}$ is a open neighborhood of $p$ contained in $D_\alpha$.
        \item It is dense: Given a basic open set $C^*$ (with $C\subseteq \omega$), let $B'\in [C]^\omega$ be such that $f_\alpha|B'$ is constant or injective. If $f_\alpha|B'$ is constant, by letting $B=B'$ we see that $B^*\subseteq C^*\cap D_\alpha$. If $f_\alpha|B'$ is injective, there exists $a\in \mathcal A$ such that $f[B']\cap a$ is infinite. Let $B\in [B']^\omega$ be such that $f[B]\subseteq a$. Then $B^*\subseteq C^*\cap D_\alpha$.
    \end{itemize}
    
    Thus fix $p\in \bigcap_{\alpha<\kappa} D_\alpha$ and $\alpha<\kappa$. We will see that $U^\alpha$ has a $p$-limit. As $p \in D_\alpha$, we have two cases:

    \textbf{Case 1:} There exists $B \in p$ such that $f_\alpha|B$ is constant. In this case, let $k \in \omega$ be such that $f[B]=\{k\}$. Then $k$ is a $p$-limit point of $U^\alpha$ as $B\subseteq \{n \in \omega: k \in U^\alpha(n)\}$ is in $p$.

    \textbf{Case 2:} There exists $B \in p$ and $a \in \mathcal A$ such that $f_\alpha|B$ is injective and $f_\alpha[B]\subseteq a$. Then $a$ is a $p$-limit point of $(U^\alpha(n): n \in \omega)$: a basic open neighborhood of $a$ is of the form $V=\{a\}\cup(a\setminus F)$, where 
    $F\subseteq a$ is finite. As $f_\alpha|B$ is injective and into $a$, there exists a finite $K\subseteq \omega$ such that $f_\alpha[B\setminus K]\subseteq a\setminus F$, thus $B\setminus K\subseteq \{n \in \omega: f_\alpha(n)\in V\}\subseteq \{n \in \omega: V\cap U_n\neq \emptyset\}$ is in $p$.
\end{proof}

\begin{corollary}The following are equivalent:

\begin{enumerate}
    \item $n(\omega^*)>\mathfrak c$,
    
    \item $\text{FA}_\mathfrak c([\omega]^\omega, \subseteq^*)$ (for every collection of $\mathfrak c$ dense subsets of the pre-ordered set $([\omega]^\omega, \subseteq^*)$ there exists a filter intersecting all of them).
    \item For every MAD family $\mathcal A$, $\Psi(\mathcal A)$ is $(\mathfrak c, \omega^*)$-pseudocompact.

    \item For every MAD family $\mathcal A$ there exists $p \in \omega^*$ such that $\Psi(\mathcal A)$ is $p$-pseudocompact.
    
    \item For every MAD family $\mathcal A$ there exists $p \in \omega^*$ such that $\exp(\Psi(\mathcal A))$ is $p$-pseudocompact.
    \item Every MAD family is pseudocompact
    \item $h=\mathfrak c$ and every base tree has a branch of cardinality $\mathfrak c$.
\end{enumerate}
\end{corollary}

\begin{proof}
    The equivalence between (1), (2) and (7) is from\cite{balcar1980space}, and a proof for the equivalence of (1), (2), (6), and (7) can be found in \cite{NewPaperPseudocompactness}. It remains to fit (3), (4) and (5) in the chain of implications.
    
    (1) implies (3) due to the previous proposition.
    
    (3) and (4) are equivalent by proposition \ref{cappseudocompact}.  
    
    (4) and (5) are equivalent by \cite[Theorem 2.4.]{ginsburg1975some}.
    
    (4) implies (6) since, in general, $p$-pseudocompact Tychonoff spaces are pseudocompact. 
\end{proof}

We do not know what is the full relation between $(\kappa, \omega^*)$-pseudocompactness, $p$-pseudocompactness and the pseudocompactness of the hyperspace. We pose the following questions:

\begin{question}
If $\exp(\Psi(\mathcal A))$ is pseudocompact, is $\Psi(\mathcal A)$ $p$-pseudocompact for some $p \in \omega^*$?
\end{question}

\begin{question}
Can we characterize the minimal $\kappa$ such that:

$$\forall \mathcal A\text{ MAD family }\Big[\left(\Psi(\mathcal A) \text{ is } (\kappa, \omega^*)-\text{pseudocompact}\implies \mathcal A \text{ is pseudocompact}\right)\Big]?$$
\end{question}

Our results show that $\kappa\leq \mathfrak c$.

\section*{Acknowledgments}
This research was developed while the authors were part of the Thematic Program on Set Theoretic Methods in Algebra, Dynamics and Geometry at the Fields Institute. The authors acknowledgment Michael Hr\v{u}\'sak and Paul Szeptycki for stimulating conversations about the content of this paper.

		\bibliographystyle{plain}
		\bibliography{bibliografia}{}

\end{document}